\documentclass[a4paper,12pt]{amsart}
\usepackage{amsfonts}
\usepackage{amssymb}
\usepackage{ifthen}
\usepackage{amscd}
\usepackage{amsxtra}
\usepackage{graphicx}
\usepackage{color}
\nonstopmode \numberwithin{equation}{section}
\newtheorem{thm}{Theorem}[section]
\newtheorem{lem}{Lemma}[section]
\newtheorem{cor}{Corollary}[section]

\newtheorem{step}{Step}[section]

\newtheorem{cl}{Claim}[section]
\newtheorem{ca}{Case}
\newtheorem{sca}{Subcase}[section]
\newtheorem{scl}{Subclaim}[section]
\newtheorem{conj}{Conjecture}

\theoremstyle{definition}
\newtheorem{defn}{Definition}[section]
\newtheorem{example}{Example}[section]
\newtheorem{op}[equation]{Open Problem}
\newtheorem{ques}[equation]{Question}
\newtheorem{rem}{Remark}[section]
\newtheorem{exam}[equation]{Example}

\newcounter {own}
\def\theown {\thesection       .\arabic{own}}

\newenvironment{pf}[1][]{%
 \vskip 3mm
 \noindent
 \ifthenelse{\equal{#1}{}}%
  {{\slshape Proof. }}%
  {{\slshape #1.} }%
 }%
{\qed\bigskip}

\newcounter{alphabet}
\newcounter{tmp}
\newenvironment{Thm}[1][]{\refstepcounter{alphabet}%
\bigskip%
\noindent%
{\bf Theorem \Alph{alphabet}}%
\ifthenelse{\equal{#1}{}}{}{ (#1)}%
{\bf .} \itshape}{\vskip 8pt}

\makeatletter
\newcommand{\Ref}[1]{\@ifundefined{r@#1}{}{\setcounter{tmp}{\ref{#1}}\Alph{tmp}}}
\makeatother

\newcommand{\IT}{{\mathbb T}}

\newcommand{\ID}{{\mathbb D}}

\def\be{\begin{equation}}
\def\ee{\end{equation}}

\newcommand{\bee}{\begin{enumerate}}
\newcommand{\eee}{\end{enumerate}}

\newcommand{\blem}{\begin{lem}}
\newcommand{\elem}{\end{lem}}
\newcommand{\bthm}{\begin{thm}}
\newcommand{\ethm}{\end{thm}}
\newcommand{\bcor}{\begin{cor}}
\newcommand{\ecor}{\end{cor}}
\newcommand{\beg}{\begin{exam}}
\newcommand{\eeg}{\end{exam}}
\newcommand{\begs}{\begin{examples}}
\newcommand{\eegs}{\end{examples}}
\newcommand{\bdefe}{\begin{defn}}
\newcommand{\edefe}{\end{defn}}
\newcommand{\bprob}{\begin{prob}}
\newcommand{\eprob}{\end{prob}}
\newcommand{\bques}{\begin{ques}}
\newcommand{\eques}{\end{ques}}
\newcommand{\bei}{\begin{itemize}}
\newcommand{\eei}{\end{itemize}}
\newcommand{\bcon}{\begin{conj}}
\newcommand{\econ}{\end{conj}}
\newcommand{\bop}{\begin{op}}
\newcommand{\eop}{\end{op}}

\newcommand{\bca}{\begin{ca}}
\newcommand{\eca}{\end{ca}}
\newcommand{\bsca}{\begin{sca}}
\newcommand{\esca}{\end{sca}}

\newcommand{\bcl}{\begin{cl}}
\newcommand{\ecl}{\end{cl}}

\newcommand{\bst}{\begin{step}}
\newcommand{\est}{\end{step}}

\newcommand{\bscl}{\begin{scl}}
\newcommand{\escl}{\end{scl}}

\newcommand{\bcons}{\begin{conjs}}
\newcommand{\econs}{\end{conjs}}
\newcommand{\bprop}{\begin{propo}}
\newcommand{\eprop}{\end{propo}}
\newcommand{\br}{\begin{rem}}
\newcommand{\er}{\end{rem}}
\newcommand{\brs}{\begin{rems}}
\newcommand{\ers}{\end{rems}}
\newcommand{\bo}{\begin{obser}}
\newcommand{\eo}{\end{obser}}
\newcommand{\bos}{\begin{obsers}}
\newcommand{\eos}{\end{obsers}}
\newcommand{\bpf}{\begin{pf}}
\newcommand{\epf}{\end{pf}}
\newcommand{\ba}{\begin{array}}
\newcommand{\ea}{\end{array}}
\newcommand{\beq}{\begin{eqnarray}}
\newcommand{\beqq}{\begin{eqnarray*}}
\newcommand{\eeq}{\end{eqnarray}}
\newcommand{\eeqq}{\end{eqnarray*}}

\newcommand{\ds}{\displaystyle}

\newcounter{minutes}
\divide\time by 60
\newcounter{hours}
\multiply\time by 60 \addtocounter{minutes}{-\time}

\begin{document}

\bibliographystyle{amsplain}
\title []
{ Bi-Lipschitz characteristic of quasiconformal self-mappings of the
unit disk satisfying bi-harmonic equation}

\def\thefootnote{}
\footnotetext{ \texttt{\tiny File:~\jobname .tex,
          printed: \number\day-\number\month-\number\year,
          \thehours.\ifnum\theminutes<10{0}\fi\theminutes}
} \makeatletter\def\thefootnote{\@arabic\c@footnote}\makeatother

\author{Shaolin Chen}
 \address{Sh. Chen, College of Mathematics and
Statistics, Hengyang Normal University, Hengyang, Hunan 421008,
People's Republic of China.} \email{mathechen@126.com}




\author{Xiantao Wang${}^{~\mathbf{*}}$}
\address{X. Wang, Department of Mathematics, Shantou University, Shantou,
Guangdong 515063, People's Republic of China.}
\email{xtwang@stu.edu.cn}

\subjclass[2000]{Primary:  30C62; Secondary: 31A05.}
 \keywords{Lipschitz continuity, bi-Lipschitz continuity, quasiconformal mapping, biharmonic equation.
 \\
${}^{\mathbf{*}}$ Corresponding author}


\begin{abstract} Suppose that $f$ is a $K$-quasiconformal self-mapping of the unit
disk $\mathbb{D}$, which satisfies the following: $(1)$ the biharmonic
equation $\Delta(\Delta f)=g$ $(g\in \mathcal{C}(\overline{\mathbb{D}}))$,
(2) the boundary condition $\Delta f=\varphi$
($\varphi\in\mathcal{C}(\mathbb{T})$ and $\mathbb{T}$ denotes the
unit circle), and $(3)$ $f(0)=0$. The purpose of this paper is to
prove that $f$ is Lipschitz continuos, and, further, it is
bi-Lipschitz continuous when $\|g\|_{\infty}$ and
$\|\varphi\|_{\infty}$ are small enough. Moreover, the estimates are
asymptotically sharp as $K\to 1$, $\|g\|_{\infty}\to 0$ and
$\|\varphi\|_{\infty}\to 0$, and thus, such a mapping $f$ behaves
almost like a rotation for sufficiently small $K$, $\|g\|_{\infty}$
and $\|\varphi\|_{\infty}$.
\end{abstract}


\maketitle \pagestyle{myheadings} \markboth{ Shaolin Chen and Xiantao Wang}{
Bi-Lipschitz characteristic of quasiconformal self-mappings of the
unit disk}

\section{Preliminaries and  main results }\label{csw-sec1}
Let $\mathbb{C}  \cong \mathbb{R}^{2}$ be the complex plane. For
$a\in\mathbb{C}$ and  $r>0$, let $\ID(a,r)=\{z:\, |z-a|<r\}$, the open disk with center $a$ and radius $r$. For convenience,
we use $\mathbb{D}_r$ to denote $\mathbb{D}(0,r)$ and $\mathbb{D}$ the open unit disk $\ID_1$. Let
$\IT$ be the unit circle, i.e., the boundary $\partial\mathbb{D}$ of $\mathbb{D}$ and $\overline{\ID}=\ID\cup \IT$. Also, we
denote by $\mathcal{C}^{m}(D)$ the set of all complex-valued
$m$-times continuously differentiable functions from $D$ into
$\mathbb{C}$, where $D$ is a subset of $\mathbb{C}$  and
$m\in\mathbb{N}_0:=\mathbb{N}\cup\{0\}$. In particular, let
$\mathcal{C}(D):=\mathcal{C}^{0}(D)$, the set of all continuous
functions in $D$.

For a real $2\times2$ matrix $A$, we use the matrix norm
$$\|A\|=\sup\{|Az|:\,|z|=1\}$$ and the matrix function
$$\lambda(A)=\inf\{|Az|:\,|z|=1\}.$$

For $z=x+iy\in\mathbb{C}$, the
formal derivative of a complex-valued function $f=u+iv$ is given
by
$$D_{f}=\left(\begin{array}{cccc}
\ds u_{x}\;~~ u_{y}\\[2mm]
\ds v_{x}\;~~ v_{y}
\end{array}\right).
$$
Then,
$$\|D_{f}\|=|f_{z}|+|f_{\overline{z}}| ~\mbox{ and }~ \lambda(D_{f})=\big| |f_{z}|-|f_{\overline{z}}|\big |,
$$
where $$f_{z}=\frac{\partial f}{\partial z}=\frac{1}{2}\big(
f_x-if_y\big)~\mbox{and}~ f_{\overline{z}}=\frac{\partial
f}{\partial \overline{z}}=\frac{1}{2}\big(f_x+if_y\big).$$

Moreover, we use
$$J_{f}:=\det D_{f} =|f_{z}|^{2}-|f_{\overline{z}}|^{2}
$$
to denote the {\it Jacobian} of $f$.

For $z, \zeta\in\mathbb{D}$ with $z\neq \zeta$, let
$$G(z,\zeta)=\log\left|\frac{1-z\overline{\zeta}}{z-\zeta}\right|~\mbox{ and}
~P(z,e^{it})=\frac{1-|z|^{2}}{|1-ze^{-it}|^{2}}$$ be the {\it
 Green function} and  the {\it Poisson kernel},
respectively, where $t\in[0,2\pi].$

Let $g\in L^{1}(\mathbb{D})$ and $f\in\mathcal{C}^{4}(\mathbb{D})$.
Of particular interest for our investigation is the following {\it
bi-harmonic equation}:

\be\label{eq-ch-1} \Delta(\Delta f)=g~\mbox{in}~\mathbb{D}\ee
with the following associated {\it Dirichlet boundary value condition}:
\be\label{eq-ch-2}
\begin{cases}
\displaystyle \Delta f=\varphi &\mbox{ in }\, \mathbb{T},\\
\displaystyle f=f^{\ast}&\mbox{ in }\, \mathbb{T},
\end{cases}\ee where $f^{\ast},~\varphi\in \mathcal{C}(\mathbb{T})$ and $$\Delta
f:=\frac{\partial^{2}f}{\partial
x^{2}}+\frac{\partial^{2}f}{\partial y^{2}}=4f_{z \overline{z}}$$
stands for the {\it Laplacian} of $f$.

By \cite[Theorem 1]{HK} (or \cite[Theorem 1]{Be}), we see that all
solutions to the equation (\ref{eq-ch-1}) satisfying the
condition (\ref{eq-ch-2}) are given by

\be\label{eq-chen1}f(z)=\mathcal{P}_{f^{\ast}}(z)+G_{1}[\varphi](z)-
G_{2}[g](z),\ee where
$$\mathcal{P}_{f^{\ast}}(z)=\frac{1}{2\pi}\int_{0}^{2\pi}P(z,e^{i\theta})f^{\ast}(e^{i\theta})d\theta,$$
 \be\label{wed-2}
 G_{1}[\varphi](z)=\frac{1}{8\pi}\int_{0}^{2\pi}(1-|z|^{2})
\left[1+\frac{\log(1-ze^{-i\theta})}{ze^{-i\theta}}+\frac{\log(1-\overline{z}e^{i\theta})}{\overline{z}e^{i\theta}}\right]\varphi(e^{i\theta})d\theta,
\ee
\beq\label{wed-3}
G_{2}[g](z)&=&\frac{1}{16\pi}\int_{\mathbb{D}}\bigg\{2|\zeta-z|^{2}G(z,\zeta)+(1-|z|^{2})(1-|\zeta|^{2})\\ \nonumber
 &&
\times\left[\frac{\log(1-z\overline{\zeta})}{z\overline{\zeta}}+\frac{\log(1-\overline{z}\zeta)}{\overline{z}\zeta}\right]\bigg\}g(\zeta)d\sigma(\zeta),
\eeq and $d\sigma$ denotes the Lebesgue area measure in
$\mathbb{D}$.  We refer the reader to \cite{HH-2, HH-1, May} etc for
more discussions in this line.

Given a subset $\Omega$ of $\mathbb{C}$, a function
$\psi:~\Omega\rightarrow\mathbb{C}$ is said to belong to the {\it
Lipschitz space} $\Lambda(\Omega)$ if
$$\sup_{z_{1},z_{2}\in\Omega,z_{1}\neq z_{2}}\frac{|\psi(z_{1})-\psi(z_{2})|}{|z_{1}-z_{2}|}<\infty.$$
Further, a function $\psi\in\Lambda(\Omega)$ is said to be {\it
bi-Lipschitz continous} if there is a positive constant $M$ such that for all
$z_{1},~z_{2}\in\Omega$, \be\label{mon-1}
M|z_{1}-z_{2}|\leq|\psi(z_{1})-\psi(z_{2})|.\ee 


For a given domain $\Omega$, we say that a function
$u:~\Omega\rightarrow\mathbb{R}$ is {\it absolutely continuous on
lines}, $ACL$ in brief, if for every closed rectangle
$R\subset\Omega$ with sides parallel to the axes $x$ and $y$,
respectively, $u$ is absolutely continuous on almost every
horizontal line and almost every vertical line in $R$. It is
well-known that partial derivatives $u_{x}$ and $u_{y}$ exist almost
everywhere in $\Omega$.

The definition carries over to complex-valued functions.

\begin{defn} Let $K\geq1$ be a constant. A sense-preserving homeomorphism $f:$ $\Omega\rightarrow D$, between
domains $\Omega$ and $D$ in $\mathbb{C}$, is {\it
$K$-quasiconformal}
 if $f$ is $ACL$ in $\Omega$ and
$$|f_{\overline{z}}|\leq\frac{K-1}{K+1}|f_{z}|\;\;\big(\mbox{or}\;\;K^{-1}\|D_{f}\|^{2}\leq J_{f}\leq K\lambda^{2}(D_{f})\big)$$
 almost everywhere in
$\Omega$.
\end{defn}

The following is the so-called Mori's Theorem (cf. \cite{CL, FV,K-M, K2, Mo}).

\begin{Thm}\label{Mori}
Suppose that $f$ is a $K$-quasiconformal self-mapping of $\mathbb{D}$ with
$f(0)=0$. Then, there exists a constant $Q(K)$, satisfying
the condition $Q(K)\rightarrow1$ as $K\rightarrow1$, such
that
$$|f(z_{2})-f(z_{1})|\leq Q(K)|z_{2}-z_{1}|^{\frac{1}{K}},$$
where the notation $Q(K)$ means that the constant $Q$ depends only on $K$.
\end{Thm}
We remark that in \cite{Qiu} it is proved $$1\leq
Q(K)\leq16^{1-\frac{1}{K}}\min\left\{\bigg(\frac{23}{8}\bigg)^{1-\frac{1}{K}},~\big(1+2^{3-2K}\big)^{\frac{1}{K}}\right\}.$$

A natural problem is that under which condition(s) a quasiconformal mapping is Lipschitz continuous.
Recently, the study of this problem has been attracted much attention.
For example, the Lipschitz characteristic
of harmonic quasiconformal mappings has been discussed by many authors (\cite{CPW-211, K4, K5, K7, KMa, MV1, PS, Pav2}).
The Lipschitz continuity of $(K, K')$-quasiconformal
harmonic mappings has also been investigated in \cite{CLSW,K8}. See, e.g., \cite{GI-1, GI-2, H-S, K3, K9, Ma, Pav2, TW, TW-1}
for more discussions on the properties of harmonic
quasiconformal mappings.
On the study of the Lipschitz continuity of quasiconformal
mappings satisfying certain elliptic PDEs, we refer to
\cite{ATM,K6,K-M,K2}. The following result is from \cite{K2}, which is
a generalization of
the main results of \cite{PS,Pav2}.

\begin{Thm}{\rm (\cite[Theorem 1.2]{K2})}\label{K-P} Suppose that $K\geq1$ is arbitrary and
$g\in\mathcal{C}(\overline{\mathbb{D}})$. Then, there exist constants
$N(K,g)$ and $M(K)$ with $\lim_{K\rightarrow1}M(K)=1$  such that if
$f$ is a $K$-quasiconformal self-mapping of $\mathbb{D}$ satisfying
the $PDE:$ $\Delta f=g$ with $f(0)=0$, then for
$z_{1},~z_{2}\in\mathbb{D}$,
$$\left(\frac{1}{M(K)}-\frac{7\|g\|_{\infty}}{6}\right)|z_{1}-z_{2}|\leq|f(z_{1})-f(z_{2})|\leq \left(M(K)+N(K,g)\|g\|_{\infty}\right)|z_{1}-z_{2}|,$$
where $\|g\|_{\infty}=\sup_{z\in\mathbb{D}}\{|g(z)|\}$.
\end{Thm}

The aim of this paper is to discuss the Lipschitz continuity of
quasiconformal self-mapping of $\mathbb{D}$ satisfying the equation \eqref{eq-ch-1}
with the boundary condition \eqref{eq-ch-2}. Our result is as follows.

\begin{thm} \label{thm-1.1}
Let $g\in\mathcal{C}(\overline{\mathbb{D}})$, $\varphi\in\mathcal{C}(\mathbb{T})$,
and let $K\geq1$ be a
constant. Suppose that $f$ is a $K$-quasiconformal self-mapping of
$\mathbb{D}$ satisfying the equation \eqref{eq-ch-1} with  $\Delta
f=\varphi$ in $\mathbb{T}$ and $f(0)=0$. Then, there are nonnegative constants
$M_{j}(K)$ and $N_{j}(K,\varphi,g)$ {\rm($j\in\{1,2\}$)} with
$$\lim_{K\rightarrow1}M_{j}(K)=1~\mbox{and}~\lim_{\|\varphi\|\rightarrow0,\|g\|\rightarrow0}N_{j}(K,\varphi,g)=0$$
such that for all $z_{1}$ and $z_{2}$ in $\ID$,
$$
\big(M_{1}(K)-N_{1}(K,\varphi,g)\big)|z_{1}-z_{2}|\leq|f(z_{1})-f(z_{2})|\leq
\big(M_{2}(K)+N_{2}(K,\varphi,g)\big)|z_{1}-z_{2}|.$$
\end{thm}

\begin{rem}
By the discussions in Step 3 of the proof of Theorem \ref{thm-1.1} in Section \ref{sec-3}, we see that the co-Lipschitz continuity coefficient
$M_{1}(K)-N_{1}(K,\varphi,g)$ is
positive for small enough norms
$\|g\|_{\infty}$ and $\|\varphi\|_{\infty}$ (for example, if $\|g\|_{\infty}\leq \frac{60}{(25+61K^{2})46^{2(K-1)}}$ and $\|\varphi\|_{\infty}\leq \frac{25}{(38+101K^{2})46^{2(K-1)}}$ (see Corollary \ref{tue-1})).
Example \ref{wed-5} shows that this condition for $f$ to be co-Lipschtz continuous cannot be replaced by the one that
$\varphi$ and $g$ are arbitrary. In Section \ref{sec-4}, another example is constructed to illustrate the possibility of $f$ from Theorem \ref{thm-1.1} to be bi-Lipschitz
continuous.
\end{rem}

We will prove several auxiliary results in the next section, Section
\ref{sec-2}. The proof of Theorem \ref{thm-1.1} will be presented in
Section \ref{sec-3}, and in Section \ref{sec-4}, two examples are
constructed.

\section{Preliminaries}\label{sec-2}
In this section, we shall prove several lemmas which will be used later on. The first lemma is as follows.

\begin{lem}\label{lem-1}
Suppose that $\varphi\in\mathcal{C}(\mathbb{T})$ and $G_{1}[\varphi]$ is defined in \eqref{wed-2}. Then, the following statements hold:

\noindent $(1)$ For $z\in\mathbb{D}$,
$$
\max\Bigg\{\left|\frac{\partial}{\partial
z}G_{1}[\varphi](z)\right|,\;\;\left|\frac{\partial}{\partial
\overline{z}}G_{1}[\varphi](z)\right|\Bigg\}
\leq\frac{\|\varphi\|_{\infty}}{4}\left[\max_{x\in[0,1)}\{h(x)\}+
\left(\frac{\pi^{2}}{3}-1\right)^{\frac{1}{2}}|z|\right],$$ where
$$
h(x)=(1-x)\left[\sum_{n=2}^{\infty}\frac{(n-1)^{2}}{n^{2}}x^{n-2}\right]^{\frac{1}{2}}\;\;\mbox{in} \;\;[0,1).$$\medskip

\noindent $(2)$ Both
$\frac{\partial}{\partial z}G_{1}[\varphi]$ and
$\frac{\partial}{\partial \overline{z}}G_{1}[\varphi]$ have
continuous extensions to the boundary, and further, for
$t\in[0,2\pi]$,
\be\label{eq-h+2}\frac{\partial}{\partial
z}G_{1}[\varphi](e^{it})=-\frac{e^{-it}}{8\pi}\int_{0}^{2\pi}
\left[1+\frac{\log(1-e^{i(t-\theta)})}{e^{i(t-\theta)}}+\frac{\log(1-e^{i(\theta-t)})}{e^{i(\theta-t)}}\right]\varphi(e^{i\theta})d\theta,
\ee
 \be\label{eq-h+3}\frac{\partial}{\partial
\overline{z}}G_{1}[\varphi](e^{it})=-\frac{e^{it}}{8\pi}\int_{0}^{2\pi}
\left[1+\frac{\log(1-e^{i(t-\theta)})}{e^{i(t-\theta)}}+\frac{\log(1-e^{i(\theta-t)})}{e^{i(\theta-t)}}\right]\varphi(e^{i\theta})d\theta,
\ee
 \be\label{eq-h+4}\left|\frac{\partial}{\partial
z}G_{1}[\varphi](e^{it})\right|\leq
\frac{\|\varphi\|_{\infty}}{4}\left(\frac{\pi^{2}}{3}-1\right)^{\frac{1}{2}},\ee
and \be\label{eq-h+5}\left|\frac{\partial}{\partial
\overline{z}}G_{1}[\varphi](e^{it})\right|\leq
\frac{\|\varphi\|_{\infty}}{4}\left(\frac{\pi^{2}}{3}-1\right)^{\frac{1}{2}}.\ee
\end{lem}

The proof of Lemma \ref{lem-1} needs the following result (cf. \cite[Proposition 2.4]{K2}).

\begin{Thm}\label{Thm-A}
Suppose that $X$ is an open subset of $\mathbb{R}$, and $(\Omega,
\mu)$ denotes a measure space. Suppose, further, that a function $F:$
$X\times\Omega\rightarrow\mathbb{R}$ satisfies the following
conditions:
\begin{enumerate}
\item[(1)] $F(x,w)$ is a measurable function of $x$ and $w$ jointly,
and is integrable with respect to $w$ for almost every $x\in X.$
\item[(2)] For almost every $w\in\Omega$, $F(x,w)$ is an absolutely
continuous function with respect to $x$. $($This guarantees that
$\partial F/\partial x$ exists almost everywhere.$)$
\item[(3)] $\partial F/\partial x$ is locally integrable, that is,
for all compact intervals $[a,b]$ contained in $X$,
$$\int_{a}^{b}\int_{\Omega}\left|\frac{\partial}{\partial x}F(x,w)\right|d\mu(w)dx<\infty.$$
\end{enumerate}
Then, $\int_{\Omega}F(x,w)d\mu(w)$ is an absolutely continuous
function with respect to $x$, and for almost every $x\in X$, its
derivative exists, which is given by
$$\frac{d}{dx}\int_{\Omega}F(x,w)d\mu(w)=\int_{\Omega}\frac{\partial}{\partial x}F(x,w)d\mu(w).$$
\end{Thm}

\subsection*{Proof of Lemma \ref{lem-1}} To prove the first statement of the lemma, we only need to show the inequality:
$$\left|\frac{\partial}{\partial
z}G_{1}[\varphi](z)\right|
\leq\frac{\|\varphi\|_{\infty}}{4}\left[\max_{x\in[0,1)}\{h(x)\}+
\left(\frac{\pi^{2}}{3}-1\right)^{\frac{1}{2}}|z|\right],$$ since the proof for the other one is similar.
Let
$$I_{1}(z)=\frac{(|z|^{2}-1)}{8\pi}\int_{0}^{2\pi}\left[\frac{1}{z(1-ze^{-i\theta})}+\frac{e^{i\theta}\log(1-ze^{-i\theta})}{z^{2}}\right]\varphi(e^{i\theta})d\theta$$
and $$I_{2}(z)=-\frac{\overline{z}}{8\pi}\int_{0}^{2\pi}
\left[1+\frac{\log(1-ze^{-i\theta})}{ze^{-i\theta}}+\frac{\log(1-\overline{z}e^{i\theta})}{\overline{z}e^{i\theta}}\right]\varphi(e^{i\theta})d\theta.$$

First, we estimate $|I_{1}(z)|$. Since
$$|I_{1}(z)|\leq \frac{(1-|z|^{2})\|\varphi\|_{\infty}}{8\pi}\int_{0}^{2\pi}\left|\sum_{n=2}^{\infty}\frac{(n-1)(ze^{-i\theta})^{n-2}e^{-i\theta}}{n}\right|d\theta,$$
and since H\"older's inequality implies
$$\frac{1}{2\pi}\int_{0}^{2\pi}\left|\sum_{n=2}^{\infty}\frac{(n-1)(ze^{-i\theta})^{n-2}e^{-i\theta}}{n}\right|d\theta\leq
 \left(\frac{1}{2\pi}\int_{0}^{2\pi}\left|\sum_{n=2}^{\infty}\frac{(n-1)(ze^{-i\theta})^{n-2}}{n}\right|^{2}d\theta\right)^{\frac{1}{2}},$$
we see that 
 $$|I_{1}(z)|\leq
\frac{(1-|z|^{2})\|\varphi\|_{\infty}}{4}\left[\sum_{n=2}^{\infty}\frac{(n-1)^{2}}{n^{2}}|z|^{2(n-2)}\right]^{\frac{1}{2}}.
$$

Let
$$h(x)=(1-x)\left[\sum_{n=2}^{\infty}\frac{(n-1)^{2}}{n^{2}}x^{n-2}\right]^{\frac{1}{2}}\;\;\mbox{in}\;\;[0,1).$$
It follows from
$$
h(x)\leq
(1-x)\left(\sum_{n=2}^{\infty}x^{n-2}\right)^{\frac{1}{2}}=(1-x)^{\frac{1}{2}}$$
that
$\max_{x\in[0,1)}\{h(x)\}$ does exist. Then, we obtain that for
$z\in\mathbb{D}$,

\be\label{kkww-1}|I_{1}(z)|\leq\frac{\|\varphi\|_{\infty}}{4}\max_{x\in[0,1)}\{h(x)\}.\ee

Next, we estimate $|I_{2}(z)|$. Since
$$|I_{2}(z)| \leq \frac{|z|\|\varphi\|_{\infty}}{8\pi}\int_{0}^{2\pi}
\left|1+\frac{\log(1-ze^{-i\theta})}{ze^{-i\theta}}+\frac{\log(1-\overline{z}e^{i\theta})}{\overline{z}e^{i\theta}}\right|d\theta,$$
we obtain from H\"older's inequality that
$$|I_{2}(z)| \leq \frac{|z|\|\varphi\|_{\infty}}{4}\left\{\frac{1}{2\pi}\int_{0}^{2\pi}
\left|1+\frac{\log(1-ze^{-i\theta})}{ze^{-i\theta}}+\frac{\log(1-\overline{z}e^{i\theta})}{\overline{z}e^{i\theta}}\right|^{2}d\theta\right\}^{\frac{1}{2}}.$$
Then, it follows from
$$\left\{\frac{1}{2\pi}\int_{0}^{2\pi}
\left|1+\frac{\log(1-ze^{-i\theta})}{ze^{-i\theta}}+\frac{\log(1-\overline{z}e^{i\theta})}{\overline{z}e^{i\theta}}\right|^{2}d\theta\right\}^{\frac{1}{2}}
\leq \left(2\sum_{n=1}^{\infty}\frac{1}{n^{2}}-1\right)^{\frac{1}{2}}$$
that

\be\label{kkww-2}|I_{2}(z)| \leq
\frac{\|\varphi\|_{\infty}}{4}\left(\frac{\pi^{2}}{3}-1\right)^{\frac{1}{2}}|z|.\ee


Now, (\ref{kkww-1}), (\ref{kkww-2}) and Theorem \Ref{Thm-A}
guarantee that

\begin{eqnarray*}
\left|\frac{\partial}{\partial
z}G_{1}[\varphi](z)\right|&=&\left|\sum_{j=1}^{2}I_{j}(z)\right|\leq
\sum_{j=1}^{2}|I_{j}(z)|\\ &\leq&
\frac{\|\varphi\|_{\infty}}{4}\max_{x\in[0,1)}\{h(x)\}+
\frac{\|\varphi\|_{\infty}}{4}\left(\frac{\pi^{2}}{3}-1\right)^{\frac{1}{2}}|z|,
\end{eqnarray*} as required.

It follows from the first statement of the lemma and the Vitali
Theorem (cf. \cite[Theorem 26.C]{Ha}) that $\frac{\partial}{\partial
z}G_{1}[\varphi]$ has a continuous extension to the boundary, and
thus,
 $$\lim_{r\rightarrow1^{-}}\frac{\partial}{\partial z}G_{1}[\varphi](re^{it})=
 \lim_{r\rightarrow1^{-}}I_{1}(re^{it})+\lim_{r\rightarrow1^{-}}I_{2}(re^{it})=I_{2}(e^{it}),$$
which implies
$$\left|\frac{\partial}{\partial
z}G_{1}[\varphi](e^{it})\right|\leq
\frac{\|\varphi\|_{\infty}}{4}\left(\frac{\pi^{2}}{3}-1\right)^{\frac{1}{2}}.$$
These show that (\ref{eq-h+2}) and
(\ref{eq-h+4}) hold.

Similarly, we see that (\ref{eq-h+3}) and (\ref{eq-h+5}) are also true.
 Hence, the proof of the lemma is complete. \qed \medskip

The following result is useful for the proof of Lemma \ref{lem-2} below.

\begin{Thm}{\rm (cf. \cite{LP})}\label{Thm-1}
For  $z\in\mathbb{D}$, we have
$$\frac{1}{2\pi}\int_{0}^{2\pi}\frac{d\theta}{|1-ze^{i\theta}|^{2\alpha}}=\sum_{n=0}^{\infty}\left(\frac{\Gamma(n+\alpha)}{n!\Gamma(\alpha)}\right)^{2}|z|^{2n},$$
where $\alpha>0$ and $\Gamma$ denotes the Gamma function.
\end{Thm}

\begin{lem}\label{lem-2} Suppose $g\in\mathcal{C}(\overline{\mathbb{D}})$ and $G_{2}[g]$ is defined in \eqref{wed-3}.
Then, the following statements hold:

\noindent $(1)$
 For
$z\in\mathbb{D}$, $$\max\Bigg\{\left|\frac{\partial}{\partial z}
G_{2}[g](z)\right|,~\left|\frac{\partial}{\partial \overline{z}}
G_{2}[g](z)\right|\Bigg\}\leq\|g\|_{\infty}\left[\frac{1}{16}+\frac{(1-|z|^{2})^{\frac{1}{2}}}{60}+\frac{2^{\frac{1}{2}}\left(1+\frac{\pi^{2}}{6}\right)^{\frac{1}{2}}}{32}|z|\right].$$

\noindent $(2)$  Both $\frac{\partial}{\partial z} G_{2}[g]$ and
$\frac{\partial}{\partial\overline{ z}} G_{2}[g]$ have continuous
extensions to the boundary, and further, for $\theta\in[0,2\pi]$,
\beq\label{eq-b1} \frac{\partial}{\partial z}
G_{2}[g](e^{i\theta})&=&\frac{1}{8\pi}\int_{\mathbb{D}}|\zeta-e^{i\theta}|^{2}\frac{\partial}{\partial
z}G(e^{i\theta},\zeta)g(\zeta)d\sigma(\zeta)\\ \nonumber
&&-\frac{e^{-i\theta}}{16\pi}\int_{\mathbb{D}}(1-|\zeta|^{2})
\bigg[\frac{\log(1-e^{i\theta}\overline{\zeta})}{e^{i\theta}\overline{\zeta}}\\
\nonumber &&
+\frac{\log(1-e^{-i\theta}\zeta)}{e^{-i\theta}\zeta}\bigg]g(\zeta)d\sigma(\zeta),
\eeq

\beq\label{eq-b2}\frac{\partial}{\partial \overline{z}}
G_{2}[g](e^{i\theta})&=&\frac{1}{8\pi}\int_{\mathbb{D}}|\zeta-e^{i\theta}|^{2}\frac{\partial}{\partial
\overline{z}}G(e^{i\theta},\zeta)g(\zeta)d\sigma(\zeta)\\ \nonumber
&&-\frac{e^{i\theta}}{16\pi}\int_{\mathbb{D}}(1-|\zeta|^{2})
\bigg[\frac{\log(1-e^{i\theta}\overline{\zeta})}{e^{i\theta}\overline{\zeta}}\\ \nonumber
&&+\frac{\log(1-e^{-i\theta}\zeta)}{e^{-i\theta}\zeta}\bigg]g(\zeta)d\sigma(\zeta),
\eeq

\be\label{eq-b3} \left|\frac{\partial}{\partial z}
G_{2}[g](e^{i\theta})\right|\leq\frac{\|g\|_{\infty}}{32}\left[1+2^{\frac{1}{2}}\left(1+\frac{\pi^{2}}{6}\right)^{\frac{1}{2}}\right],\ee
and \be\label{eq-b4} \left|\frac{\partial}{\partial \overline{z}}
G_{2}[g](e^{i\theta})\right|\leq\frac{\|g\|_{\infty}}{32}\left[1+2^{\frac{1}{2}}\left(1+\frac{\pi^{2}}{6}\right)^{\frac{1}{2}}\right].\ee
\end{lem}

\bpf To prove the first statement, we only need to prove the inequality:
$$\left|\frac{\partial}{\partial z}
G_{2}[g](z)\right|\leq\|g\|_{\infty}\left[\frac{1}{16}+\frac{(1-|z|^{2})^{\frac{1}{2}}}{60}+\frac{2^{\frac{1}{2}}\left(1+\frac{\pi^{2}}{6}\right)^{\frac{1}{2}}}{32}|z|\right]$$
 because the proof for the other one is similar. For
this, let
$$I_{3}(z)=\frac{1}{8\pi}\int_{\mathbb{D}}(\overline{z}-\overline{\zeta})G(z,\zeta)g(\zeta)d\sigma(\zeta),$$
$$I_{4}(z)=\frac{1}{8\pi}\int_{\mathbb{D}}|\zeta-z|^{2}\frac{\partial}{\partial z}G(z,\zeta)g(\zeta)d\sigma(\zeta),$$
$$I_{5}(z)=-\frac{1}{16\pi}\int_{\mathbb{D}}\overline{z}(1-|\zeta|^{2})
\left[\frac{\log(1-z\overline{\zeta})}{z\overline{\zeta}}+\frac{\log(1-\overline{z}\zeta)}{\overline{z}\zeta}\right]g(\zeta)d\sigma(\zeta),$$
and
$$I_{6}(z)=-\frac{1}{16\pi}\int_{\mathbb{D}}(1-|\zeta|^{2})(1-|z|^{2})\left[\frac{1}{z(1-z\overline{\zeta})}
+\frac{\log(1-z\overline{\zeta})}{z^{2}\overline{\zeta}}\right]g(\zeta)d\sigma(\zeta).$$

We are going to estimate the norms of $I_{3}(z)$, $I_{4}(z)$, $I_{5}(z)$, $I_{6}(z)$, respectively. Before these estimates, we need some preparation.
Set $$w=\frac{z-\zeta}{1-\overline{z}\zeta}.$$ Then,
\be\label{eq-t6}|J_{\zeta}(w)|=\frac{(1-|z|^{2})^{2}}{|1-\overline{z}w|^{4}},\;\;\zeta-z=\frac{w(|z|^{2}-1)}{1-\overline{z}w}, \ee
and
\be\label{eq-t6.02}1-|\zeta|^{2}=\frac{(1-|w|^{2})(1-|z|^{2})}{|1-\overline{z}w|^{2}}.\ee

Firstly, we estimate $|I_{3}(z)|$. Since (\ref{eq-t6}) and
(\ref{eq-t6.02}) guarantee that
$$
|I_{3}(z)|\leq \frac{\|g\|_{\infty}}{8\pi}(1-|z|^{2})^{3}\int_{\mathbb{D}}\left(|w|\log\frac{1}{|w|}\right)\frac{d\sigma(w)}{|1-\overline{z}w|^{5}},
$$ by letting $w=re^{i\theta}$, we obtain
$$
|I_{3}(z)|\leq \frac{\|g\|_{\infty}}{4}(1-|z|^{2})^{3}\int_{0}^{1}r^{2}\log\frac{1}{r}
\left(\frac{1}{2\pi}\int_{0}^{2\pi}\frac{d\theta}{|1-\overline{z}re^{i\theta}|^{5}}\right)dr.
$$

Moreover, H\"older's inequality and Theorem \Ref{Thm-1} show that
$$
\frac{1}{2\pi}\int_{0}^{2\pi}\frac{d\theta}{|1-\overline{z}re^{i\theta}|^{5}}\leq
\left(\frac{1}{2\pi}\int_{0}^{2\pi}\frac{d\theta}{|1-\overline{z}re^{i\theta}|^{4}}\right)^{\frac{1}{2}}
\left(\frac{1}{2\pi}\int_{0}^{2\pi}\frac{d\theta}{|1-\overline{z}re^{i\theta}|^{6}}\right)^{\frac{1}{2}}.
$$
Since
$$
\left(\frac{1}{2\pi}\int_{0}^{2\pi}\frac{d\theta}{|1-\overline{z}re^{i\theta}|^{4}}\right)^{\frac{1}{2}}=\left[\sum_{n=0}^{\infty}(n+1)^{2}|zr|^{2n}\right]^{\frac{1}{2}}
$$
and
$$
\left(\frac{1}{2\pi}\int_{0}^{2\pi}\frac{d\theta}{|1-\overline{z}re^{i\theta}|^{6}}\right)^{\frac{1}{2}}
=\left[\sum_{n=0}^{\infty}(n+1)^{2}\Big(\frac{n}{2}+1\Big)^{2}|zr|^{2n}\right]^{\frac{1}{2}},
$$
we see that
$$
\frac{1}{2\pi}\int_{0}^{2\pi}\frac{d\theta}{|1-\overline{z}re^{i\theta}|^{5}}
\leq\frac{1}{4}\sum_{n=0}^{\infty}(n+1)^{2}(n+2)^{2}|zr|^{2n},
$$
which implies
\be\label{eq-t9} |I_{3}(z)|\leq
\frac{\|g\|_{\infty}}{64}(1-|z|^{2})^{3}\sum_{n=0}^{\infty}(n+1)(n+2)|z|^{2n}=\frac{\|g\|_{\infty}}{32}.
\ee

Secondly, we estimate $|I_{4}(z)|$. By Theorem \Ref{Thm-1}, we obtain that

$$\frac{1}{2\pi}
\int_{0}^{2\pi}\frac{1}{|1-\overline{z}re^{i\theta}|^{6}}d\theta=\sum_{n=0}^{\infty}\frac{(n+1)^{2}(n+2)^{2}}{4}|z|^{2n}r^{2n},$$
which implies

 \be\label{wwkk-3} \int_{0}^{1}\left(\frac{1}{2\pi}
\int_{0}^{2\pi}\frac{r^{2}(1-r^{2})}{|1-\overline{z}re^{i\theta}|^{6}}d\theta\right)dr
\leq \frac{1}{4(1-|z|^{2})^{3}}.\ee

Moreover, by (\ref{eq-t6}), (\ref{eq-t6.02}), and by applying $w=re^{i\theta}$, we have
\begin{eqnarray*}
\frac{1}{2\pi}\int_{\mathbb{D}}\frac{|z-\zeta|(1-|\zeta|^{2})}{|1-z\overline{\zeta}|}d\sigma(\zeta)
=(1-|z|^{2})^{3}\int_{0}^{1}\left(\frac{1}{2\pi}
\int_{0}^{2\pi}\frac{r^{2}(1-r^{2})}{|1-\overline{z}re^{i\theta}|^{6}}d\theta\right)dr,
\end{eqnarray*}
which, together with (\ref{wwkk-3}), yields \beq\label{eq-t10}
 |I_{4}(z)|\leq \frac{\|g\|_{\infty}}{16\pi}\int_{\mathbb{D}}\frac{|z-\zeta|(1-|\zeta|^{2})}{|1-z\overline{\zeta}|}d\sigma(\zeta)
=\frac{\|g\|_{\infty}}{32}. \eeq

Next, we estimate $|I_{5}(z)|$. Let
$$\mathcal{A}_{1}(z,\rho)=\frac{1}{2\pi}\int_{0}^{2\pi}\left|\sum_{n=1}^{\infty}\frac{(z\rho
e^{-it})^{n-1}}{n}+\sum_{n=1}^{\infty}\frac{(\overline{z}\rho
e^{it})^{n-1}}{n}\right|dt.$$

It follows from H\"older's inequality that

$$
\mathcal{A}_{1}(z,\rho)\leq \left(\frac{1}{2\pi}\int_{0}^{2\pi}\left|\sum_{n=1}^{\infty}\frac{(z\rho
e^{-it})^{n-1}}{n}+\sum_{n=1}^{\infty}\frac{(\overline{z}\rho
e^{it})^{n-1}}{n}\right|^{2}dt\right)^{\frac{1}{2}}
\leq 2^{\frac{1}{2}}\left(1+\frac{\pi^{2}}{6}\right)^{\frac{1}{2}}.
$$

By letting $\zeta=\rho e^{it}$, we get

\beq\label{eq-t11}
|I_{5}(z)|&\leq& \frac{|z|\|g\|_{\infty}}{8}\int_{0}^{1}\rho(1-\rho^{2})\mathcal{A}_{1}(z,\rho)d\rho
\leq \frac{2^{\frac{1}{2}}\left(1+\frac{\pi^{2}}{6}\right)^{\frac{1}{2}}\|g\|_{\infty}}{32}|z|.
\eeq

Finally, we estimate $|I_{6}(z)|$.
Let
$$\mathcal{A}_{2}(z)=\int_{0}^{1}\rho(1-\rho^{2})\left[\frac{1}{2\pi}\int_{0}^{2\pi}
\left|\sum_{n=2}^{\infty}\frac{(n-1)}{n}z^{n-2}\rho^{n-1}e^{-it(n-1)}\right|dt\right]d\rho.$$

By using H\"older's inequality, we obtain that
$$
\mathcal{A}_{2}(z)\leq\int_{0}^{1}\rho(1-\rho^{2})\left[\frac{1}{2\pi}\int_{0}^{2\pi}
\left|\sum_{n=2}^{\infty}\frac{(n-1)}{n}z^{n-2}\rho^{n-1}e^{-it(n-1)}\right|^{2}dt\right]^{\frac{1}{2}}d\rho.
$$

Since

\begin{eqnarray*}
\frac{1}{2\pi}\int_{0}^{2\pi}
\left|\sum_{n=2}^{\infty}\frac{(n-1)}{n}z^{n-2}\rho^{n-1}e^{-it(n-1)}\right|^{2}dt&=&
\rho^2
\sum_{n=2}^{\infty}\frac{(n-1)^{2}}{n^{2}}(|z|\rho)^{2(n-2)}\\
&\leq&\rho^2 \sum_{n=2}^{\infty}|z|^{2(n-2)},
\end{eqnarray*} we see that
$$\mathcal{A}_{2}(z) \leq \int_{0}^{1}\rho^{2}(1-\rho^{2})\left(\sum_{n=0}^{\infty}|z|^{2n}\right)^{\frac{1}{2}}d\rho=\frac{2}{15(1-|z|^{2})^{\frac{1}{2}}},$$
which imples
 \beq\label{eq-t12}
|I_{6}(z)| \leq \frac{(1-|z|^{2})\|g\|_{\infty}}{8}\mathcal{A}_{2}(z)\leq \frac{\|g\|_{\infty}(1-|z|^{2})^{\frac{1}{2}}}{60}.
 \eeq
 Therefore, by  (\ref{eq-t9}), (\ref{eq-t10}), (\ref{eq-t11}), (\ref{eq-t12}) and Theorem \Ref{Thm-A}, we conclude that

\beqq\label{eq-t14}\left|\frac{\partial}{\partial z}
G_{2}[g](z)\right|\leq\sum_{j=3}^{6}|I_{j}(z)|\leq \|g\|_{\infty}\left[\frac{1}{16}+\frac{(1-|z|^{2})^{\frac{1}{2}}}{60}+\frac{2^{\frac{1}{2}}\left(1+\frac{\pi^{2}}{6}\right)^{\frac{1}{2}}}{32}|z|\right],\eeqq
as required.

It follows from the first statement of the lemma, along with the Vitali
Theorem (cf. \cite[Theorem 26.C]{Ha}), that $\frac{\partial}{\partial
z} G_{2}[g]$ has a continuous extension to the boundary.

Since for
$\zeta\in\mathbb{D}$, we have
$$\lim_{|z|\rightarrow1^{-}}G(z,\zeta)=\lim_{|z|\rightarrow1^{-}}\log\left|\frac{1-z\overline{\zeta}}{z-\zeta}\right|=0,$$
which gives
\beqq\label{ww-1}\lim_{|z|\rightarrow1^{-}}I_{3}(z)=0,\eeqq and because
(\ref{eq-t12}) leads to
\beqq\label{ww-2}\lim_{|z|\rightarrow1^{-}}I_{6}(z)=0,\eeqq
we have
$$
 \lim_{r\rightarrow1^{-}}\frac{\partial}{\partial
z}G_{2}[g](re^{i\theta})=I_{4}(e^{i\theta})+I_{5}(e^{i\theta}).$$
Then, (\ref{eq-b3}) easily follows from (\ref{eq-t10}) and
(\ref{eq-t11}).

Similarly, we know that
(\ref{eq-b2}) and (\ref{eq-b4}) are also true. Hence, the
lemma is proved.
 \epf

\begin{lem}\label{lem-main}
For  $\varphi\in\mathcal{C}(\mathbb{T})$ and
$g\in\mathcal{C}(\overline{\mathbb{D}})$, suppose that $f$ is a
sense-preserving
 homeomorphism from $\overline{\mathbb{D}}$ onto  itself
satisfying {\rm (\ref{eq-ch-1})} and $\Delta f=\varphi$ in
$\mathbb{T}$, and suppose that $f$ is Lipschitz continuous in
$\mathbb{D}$. Then, for almost every $e^{i\theta}\in\mathbb{T}$, the following limits exist:
\be
\label{qw-1} D_{f}(e^{i\theta}):=\lim_{z\rightarrow
e^{i\theta},z\in\mathbb{D}}D_{f}(z)\;\;\mbox{and}\;\; J_{f}(e^{i\theta}):=\lim_{z\rightarrow e^{i\theta},z\in\mathbb{D}}J_{f}(z).\ee
Further, we have
\beq\label{eq-sh-1}J_{f}(e^{i\theta})&\leq&\frac{\eta'(\theta)}{2\pi}\int_{0}^{2\pi}\frac{|f(e^{it})-f(e^{i\theta})|^{2}}{|e^{it}-e^{i\theta}|^{2}}dt
+\frac{\eta'(\theta)\|\varphi\|_{\infty}}{2}\left(\frac{\pi^{2}}{3}-1\right)^{\frac{1}{2}}\\
\nonumber&&+\frac{\eta'(\theta)\|g\|_{\infty}}{16}\left[1+2^{\frac{1}{2}}\left(1+\frac{\pi^{2}}{6}\right)^{\frac{1}{2}}\right]
\eeq and
\beq\label{eq-sh-1.1}J_{f}(e^{i\theta})&\geq&\frac{\eta'(\theta)}{2\pi}\int_{0}^{2\pi}\frac{|f(e^{it})-f(e^{i\theta})|^{2}}{|e^{it}-e^{i\theta}|^{2}}dt
-\frac{\eta'(\theta)\|\varphi\|_{\infty}}{2}\left(\frac{\pi^{2}}{3}-1\right)^{\frac{1}{2}}\\
\nonumber&&-\frac{\eta'(\theta)\|g\|_{\infty}}{16}\left[1+2^{\frac{1}{2}}\left(1+\frac{\pi^{2}}{6}\right)^{\frac{1}{2}}\right],
\eeq where $f(e^{i\theta})=e^{i\eta(\theta)}$ and $\eta(\theta)$ is
a real-valued function in $[0,2\pi]$.
\end{lem}


Before the proof of Lemma \ref{lem-main}, let us recall the following result (cf. \cite[Lemma 2.1]{K2}).

\begin{Thm}\label{Lem-Ka-1}
Suppose that $f$ is a harmonic mapping defined in $\mathbb{D}$ and its
formal derivative $D_{f}$ is bounded in $\mathbb{D}$ {\rm(}or
equivalently, according to Rademacher's theorem, suppose that $f$ itself
is Lipschitz continuous in $\mathbb{D}${\rm)}. Then, there exists a
mapping $A\in L^{\infty}(\mathbb{T})$ such that
$D_{f}(z)=\mathcal{P}_{A}(z)$ and for almost every
$e^{i\theta}\in\mathbb{T}$,
$$\lim_{r\rightarrow1^{-}}D_{f}(re^{i\theta})=A(e^{i\theta}).$$

Moreover, the function $F(e^{i\theta}):=f(e^{i\theta})$ is
differentiable almost everywhere in $[0,2\pi]$ and
$$A(e^{i\theta})ie^{i\theta}=\frac{\partial}{\partial\theta}F(e^{i\theta}).$$
\end{Thm}

\subsection*{Proof of Lemma \ref{lem-main}} We first prove the existence of the two limits in (\ref{qw-1}). By Lemmas \ref{lem-1} and
\ref{lem-2}, we get that for any $e^{i\theta}\in\mathbb{D}$,

\be\label{eq-kk-1}\lim_{z\rightarrow
e^{i\theta},z\in\mathbb{D}}D_{G_{1}[\varphi]}(z)=D_{G_{1}[\varphi]}(e^{i\theta})\;\;\mbox{and}\;\;\lim_{z\rightarrow
e^{i\theta},z\in\mathbb{D}}D_{G_{2}[g]}(z)=D_{G_{2}[g]}(e^{i\theta}).\ee

Again, by Lemmas  \ref{lem-1} and \ref{lem-2}, we know that
$$\|D_{G_{1}[\varphi]}\|<\infty\;\; \mbox{and}\;\; \|D_{G_{2}[g]}\|<\infty,$$ which
implies the Lipschitz continuity of $G_{1}[\varphi]$ and $G_{2}[g]$
in $\mathbb{D}$. Since $f$ is Lipschitz continuous in $\mathbb{D}$,
we see that $\|D_{f}\|$ is bounded in $\mathbb{D}$. Thus, it follows
from (\ref{eq-chen1}) that $\mathcal{P}_{f^{\ast}}$ is also
Lipschitz continuous in $\mathbb{D}$, where
$f^{\ast}=f|_{\mathbb{T}}$. Now, we conclude from Theorem
\Ref{Lem-Ka-1} that
 for almost every $e^{i\theta}\in\mathbb{T}$,
$$\lim_{z\rightarrow
e^{i\theta},z\in\mathbb{D}}D_{\mathcal{P}_{f^{\ast}}}(z)$$ does
exist, which, together with (\ref{eq-chen1}) and (\ref{eq-kk-1}),
guarantees that for almost every $\theta\in[0,2\pi],$
$$\lim_{z\rightarrow
e^{i\theta},z\in\mathbb{D}}D_{f}(z)$$ also exists.

Since $$J_{f}(z)=\det
D_{f}(z),$$ obviously, we see that  $$\lim_{z\rightarrow
e^{i\theta},z\in\mathbb{D}}J_{f}(z)$$ exists for almost every
$\theta\in[0,2\pi].$

Next, we demonstrate the estimates in (\ref{eq-sh-1}) and  (\ref{eq-sh-1.1}). For convenience, in the rest of the proof of the lemma, let
$$D_{f}(e^{i\theta})=\lim_{z\rightarrow e^{i\theta},z\in\mathbb{D}}D_{f}(z)\;\;\mbox{and}\;\; J_{f}(e^{i\theta})=\lim_{z\rightarrow
e^{i\theta},z\in\mathbb{D}}J_{f}(z).$$


By Lebesgue Dominated
Convergence Theorem, the boundedness of $\|D_{f}\|$, and by letting
$z=re^{i\theta}\in\mathbb{D}$, we see that for any fixed
$\theta_{1}\in[0,2\pi]$, \beq\label{eq-d1}
f(e^{i\theta})&=&\lim_{r\rightarrow1^{-}}f(re^{i\theta})=\lim_{r\rightarrow1^{-}}\int_{\theta_{1}}^{\theta}\frac{\partial}{\partial
t}f(re^{it})dt+f(e^{i\theta_{1}})\\ \nonumber
&=&\int_{\theta_{1}}^{\theta}\lim_{r\rightarrow1^{-}}\left[ir\big(f_{z}(re^{it})e^{it}-f_{\overline{z}}(re^{it})e^{-it}\big)\right]dt+f(e^{i\theta_{1}}),
\eeq which implies that $f(e^{i\theta})$ is absolutely continuous.
Let $\eta(\theta)$ be a real-valued function in $[0,2\pi]$ such that
$$e^{i\eta(\theta)}=f(e^{i\theta}).$$
Then,
\be\label{eq-r1}f'(e^{i\theta})=i\eta'(\theta)e^{i\eta(\theta)}\ee
holds almost everywhere in $[0,2\pi].$

Since
$$J_{f}(re^{i\theta})=|f_{z}(re^{i\theta})|^{2}-|f_{\overline{z}}(re^{i\theta})|^{2}=
-\mbox{Re}\left(\overline{\frac{\partial f}{\partial
r}}\frac{i}{r}\frac{\partial f}{\partial \theta}\right),$$ we infer
from (\ref{eq-r1}) that

\beq\label{sh-2}
J_{f}(e^{i\theta})&=&\lim_{r\rightarrow1^{-}}J_{f}(re^{i\theta})=-\lim_{r\rightarrow1^{-}}\mbox{Re}\left(\overline{\frac{\partial
f}{\partial r}}\frac{i}{r}\frac{\partial f}{\partial
\theta}\right)=I_{7}- I_{8}+I_{9},\eeq where
$$I_{7}=\lim_{r\rightarrow1^{-}}\mbox{Re}\left(\frac{\overline{f(e^{i\theta})-\mathcal{P}_{f^{\ast}}(re^{i\theta})}}{1-r}\cdot\eta'(\theta)f(e^{i\theta})\right),$$
$$I_{8}=\lim_{r\rightarrow1^{-}}\mbox{Re}\left(\frac{\overline{G_{1}[\varphi](re^{i\theta})}}{1-r}\cdot\eta'(\theta)f(e^{i\theta})\right),$$
and
$$I_{9}=\lim_{r\rightarrow1^{-}}\mbox{Re}\left(\frac{\overline{G_{2}[g](re^{i\theta})}}{1-r}\cdot\eta'(\theta)f(e^{i\theta})\right).$$

Now, we are going to prove (\ref{eq-sh-1}) and (\ref{eq-sh-1.1}) by estimating the quantities $I_{7}$, $|I_{8}|$ and $|I_{9}|,$ respectively. We start with the estimate of $I_{7}$. Since
$$
\mbox{Re}\langle
f(e^{i\theta}),f(e^{i\theta})-f(e^{it})\rangle=\mbox{Re}\big[f(e^{i\theta})(\overline{f(e^{i\theta})-f(e^{it})})\big]
=\frac{1}{2}|f(e^{it})-f(e^{i\theta})|^{2}
 $$ and
$$I_{7}=\lim_{r\rightarrow1^{-}}\mbox{Re}\left(\frac{1}{2\pi}\int_{0}^{2\pi}\frac{1+r}{|1-re^{i(\theta-t)}|^{2}}\langle
\eta'(\theta)f(e^{i\theta}),f(e^{i\theta})-f(e^{it})\rangle\,dt\right),$$ where $\langle\cdot,\cdot\rangle$ denotes the inner product, it follows that
\be\label{eq-s3}
I_{7}=\eta'(\theta)\frac{1}{2\pi}\int_{0}^{2\pi}\frac{|f(e^{it})-f(e^{i\theta})|^{2}}{|e^{it}-e^{i\theta}|^{2}}dt.
\ee

Next, we estimate $|I_{8}|$. For this, let
$$\mathcal{A}_{3}(r)=\frac{1}{2\pi}\int_{0}^{2\pi}\bigg|1+\frac{\log(1-re^{i(\theta-t)})}{re^{i(\theta-t)}}
+\frac{\log(1-re^{i(t-\theta)})}{re^{i(t-\theta)}}\bigg|dt.$$ By
H\"older's inequality, we have

\begin{eqnarray*}
\lim_{r\rightarrow1^{-}}\mathcal{A}_{3}(r)&\leq&\lim_{r\rightarrow1^{-}}\Bigg(\frac{1}{2\pi}\int_{0}^{2\pi}
\bigg|1-\sum_{n=1}^{\infty}\frac{(re^{i(t-\theta)})^{n-1}}{n}-\sum_{n=1}^{\infty}\frac{(re^{i(\theta-t)})^{n-1}}{n}\bigg|^{2}dt\Bigg)^{\frac{1}{2}}\\
&=&\left(\frac{\pi^{2}}{3}-1\right)^{\frac{1}{2}},
\end{eqnarray*}
which yields
 \beqq
|I_{8}|&=&\bigg|\lim_{r\rightarrow1^{-}}\mbox{Re}\bigg[\frac{1}{8\pi}\int_{0}^{2\pi}\left(1+\frac{\log(1-re^{i(\theta-t)})}{re^{i(\theta-t)}}
+\frac{\log(1-re^{i(t-\theta)})}{re^{i(t-\theta)}}\right)\\
\nonumber&&\times(1+r)\langle\eta'(\theta)f(e^{i\theta}),\varphi(e^{it})\rangle\,dt\bigg]\bigg|\\
\nonumber&\leq&\frac{\eta'(\theta)\|\varphi\|_{\infty}}{2}\lim_{r\rightarrow1^{-}}\mathcal{A}_{3}(r).
\eeqq
Thus,
 \be\label{eq-s4}
|I_{8}|\leq \frac{\eta'(\theta)\|\varphi\|_{\infty}}{2}\left(\frac{\pi^{2}}{3}-1\right)^{\frac{1}{2}}.
\ee

Finally, we estimate $|I_{9}|$. To reach this goal, let
$$
\mathcal{A}_{4}= \lim_{z\rightarrow
e^{i\theta},z\in\mathbb{D}}\frac{1}{8\pi}\bigg|\int_{\mathbb{D}}|\zeta-z|^{2}\frac{(G(z,\zeta)-G(e^{i\theta},\zeta))}{1-|z|}
\langle \eta'(\theta)f(e^{i\theta}),g(\zeta)\rangle\,
d\sigma(\zeta)\bigg|
$$
and
\begin{eqnarray*}
\mathcal{A}_{5}&=&\lim_{z\rightarrow
e^{i\theta},z\in\mathbb{D}}\frac{1}{16\pi}\bigg|\int_{\mathbb{D}}(1+|z|)(1-|\zeta|^{2})\\
\nonumber&&\times
\left[\frac{\log(1-z\overline{\zeta})}{z\overline{\zeta}}+\frac{\log(1-\overline{z}\zeta)}{\overline{z}\zeta}\right]\langle
\eta'(\theta)f(e^{i\theta}),g(\zeta)\rangle\, d\sigma(\zeta)\bigg|.\end{eqnarray*}

Since $$\lim_{z\rightarrow
e^{i\theta},z\in\mathbb{D}}\frac{G(z,\zeta)}{1-|z|}=\lim_{z\rightarrow
e^{i\theta},z\in\mathbb{D}}\frac{G(z,\zeta)-G(e^{i\theta},\zeta)}{1-|z|}=P(\zeta,e^{i\theta}),$$
 we deduce that
$$
\mathcal{A}_{4}\leq\lim_{z\rightarrow
e^{i\theta},z\in\mathbb{D}}\frac{\eta'(\theta)\|g\|_{\infty}}{8\pi}\int_{\mathbb{D}}|\zeta-z|^{2}P(\zeta,e^{i\theta})d\sigma(\zeta)=\frac{\eta'(\theta)\|g\|_{\infty}}{16}
$$
and
$$
\mathcal{A}_{5}\leq \lim_{z\rightarrow
e^{i\theta},z\in\mathbb{D}}\frac{\eta'(\theta)\|g\|_{\infty}}{8\pi}\int_{\mathbb{D}}(1-|\zeta|^{2})
\left|\frac{\log(1-z\overline{\zeta})}{z\overline{\zeta}}+\frac{\log(1-\overline{z}\zeta)}{\overline{z}\zeta}\right|d\sigma(\zeta).
$$

Now, we conclude from \eqref{wed-3} and (\ref{eq-t11}) that
\be\label{eq-s5} |I_{9}|\leq\mathcal{A}_{4}+\mathcal{A}_{5}
\leq\frac{\eta'(\theta)\|g\|_{\infty}}{16}\left[1+2^{\frac{1}{2}}\left(1+\frac{\pi^{2}}{6}\right)^{\frac{1}{2}}\right].
 \ee
Hence, (\ref{eq-sh-1}) and (\ref{eq-sh-1.1}) follow from the inequalities
(\ref{sh-2}) $\sim$ (\ref{eq-s5}) along with the following chain of inequalities:
$$I_{7}-|I_{8}|-|I_{9}|\leq J_{f}(e^{i\theta})\leq I_{7}+|I_{8}|+|I_{9}|.$$
The proof of the lemma is complete. \qed
\medskip

The following is the so-called Heinz-Theorem.

\begin{Thm}  $($\cite[Theorem]{He}$)$ \label{Heinz}
Suppose that $f$ is a harmonic homeomorphism of $\mathbb{D}$ onto
itself with $f(0)=0$. Then, for $z\in\mathbb{D},$
$$|f_{z}(z)|^{2}+|f_{\overline{z}}(z)|^{2}\geq\frac{1}{\pi^{2}}.$$
\end{Thm}

Our next lemma is a generalization of Theorem \Ref{Heinz}.

\begin{lem}\label{Chen-Wang}
Suppose that $f$ is a harmonic homeomorphism of $\mathbb{D}$ onto
itself with $f(a)=0$, where $a\in\mathbb{D}.$ Then, for
$z\in\mathbb{D},$
$$|f_{z}(z)|^{2}+|f_{\overline{z}}(z)|^{2}\geq \frac{(1-|a|)^{2}}{\pi^{2}(1+|a|)^{2}}.$$
\end{lem}
\bpf For $\zeta\in\mathbb{D},$ let
$$\phi(\zeta)=\frac{\zeta+a}{1+\overline{a}\zeta}.$$
Then, $$\phi'(\zeta)=\frac{1-|a|^{2}}{(1+\overline{a}\zeta)^{2}}.$$

Let
$$\mathcal{F}(\zeta)=f(\phi(\zeta)).$$ Then, $\mathcal{F}$ is
also a harmonic homeomorphism of $\mathbb{D}$ onto itself with
$\mathcal{F}(0)=0.$ By Theorem \Ref{Heinz}, we have
$$(|f_{w}(\phi(\zeta))|^{2}+|f_{\overline{w}}(\phi(\zeta))|^{2})\frac{(1-|a|^{2})^{2}}{|1+\overline{a}\zeta|^{4}}=
|\mathcal{F}_{\zeta}(\zeta)|^{2}+|\mathcal{F}_{\overline{\zeta}}(\zeta)|^{2}\geq
\frac{1}{\pi^{2}},
$$
which implies
$$
|f_{w}(\phi(\zeta))|^{2}+|f_{\overline{w}}(\phi(\zeta))|^{2}\geq\frac{1}{\pi^{2}}\frac{|1+\overline{a}\zeta|^{4}}{(1-|a|^{2})^{2}}\geq
\frac{(1-|a|)^{2}}{\pi^{2}(1+|a|)^{2}},
$$ where $w=\phi(\zeta)$. This is what we need.
 \epf

The following result is a direct consequence of Lemma
\ref{Chen-Wang}. \bcor\label{H-lem} Suppose that $f$ is a harmonic
homeomorphism of $\mathbb{D}$ onto itself. Then
$$\inf_{z\in\mathbb{D}}|f_{z}(z)|>0.$$

\ecor

\section{The proof of Theorem \ref{thm-1.1}}\label{sec-3}
The purpose of this section is to prove Theorem \ref{thm-1.1}. The
proof consists of three steps. In the first step, the Lipschitz continuity of the mappings $f$ is proved, the co-Lipschitz continuity of $f$ is demonstrated in the second step,
and in the third step, the Lipschitz and co-Lipschitz continuity coefficients obtained in the first two steps are shown to have bounds
with the forms as required in Theorem \ref{thm-1.1}.

Before the proof, let us recall a
result due to Kalaj and Mateljevi\'c, which is used in the
discussions of the first step.

\begin{Thm}{\rm (\cite[Theorem 3.4]{K-M})}\label{KM-P}
Suppose that $f$ is a quasiconformal $\mathcal{C}^{2}$ diffeomorphism from
the plane domain $\Omega$ with $\mathcal{C}^{1,\alpha}$ compact
boundary onto the  plane domain $\Omega^{\ast}$ with
$\mathcal{C}^{2,\alpha}$ compact boundary. If there exist constants
$a_{1}$ and $b_{1}$ such that
$$|\Delta f(z)|\leq a_{1}\|D_{f}(z)\|^{2}+b_{1}$$ in $\Omega$, then $f$ has
bounded partial derivatives. In particular, it is a
Lipschitz mapping in $\Omega$.
\end{Thm}

\bst\label{bst-1} Lipschitz continuity.
 \est

We start the discussions of this step with the following claim. \bcl
The limits $$\lim_{z\rightarrow
\xi\in\mathbb{T},z\in\mathbb{D}}D_{f}(z)\;\; \mbox{and}\;\;
\lim_{z\rightarrow \xi\in\mathbb{T},z\in\mathbb{D}}J_{f}(z)$$ exist
almost everywhere in $\IT$. \ecl

We are going to verify the existence of  these two limits by
applying Theorem \Ref{KM-P} and Lemma \ref{lem-main}. For this, we
need to get an upper bound of the quantity $|\Delta f(z)|$ as stated
in \eqref{eq-x1.1} below. By the formula (1.3) in \cite{K2} (see
also \cite[pp. 118-120]{Ho}), we have that for $z\in \ID$,
\beqq\label{eq-yy2}\Delta
f(z)=\mathcal{P}_{\varphi}(z)-\frac{1}{2\pi}\int_{\mathbb{D}}G(z,\zeta)g(\zeta)d\sigma(\zeta).
\eeqq

Since Theorem \Ref{Thm-1} implies
$$\frac{1}{2\pi}\int_{0}^{2\pi}\frac{dt}{|1-\overline{z}re^{it}|^{4}}=\sum_{n=0}^{\infty}(n+1)^{2}|z|^{2n}r^{2n},$$
by letting
$w=\frac{z-\zeta}{1-\overline{z}\zeta}=re^{it}$, we obtain
\beqq\label{yy-1}\frac{1}{2\pi}\int_{\mathbb{D}}G(z,\zeta)d\sigma(\zeta)=\frac{1}{2\pi}\int_{\mathbb{D}}\left(\log\frac{1}{|w|}\right)
\frac{(1-|z|^{2})^{2}}{|1-\overline{z}w|^{4}}d\sigma(w)=\frac{(1-|z|^{2})}{4},\eeqq
and so, we get
\be\label{eq-x1.1} |\Delta f(z)|\leq |\mathcal{P}_{\varphi}(z)|+\frac{\|g\|_{\infty}}{2\pi}\int_{\mathbb{D}}G(z,\zeta)d\sigma(\zeta)
\leq \|\varphi\|_{\infty}+\frac{\|g\|_{\infty}}{4}.\ee
Now, the existence of the limits
$$D_{f}(\xi)=\lim_{z\rightarrow
\xi\in\mathbb{T},z\in\mathbb{D}}D_{f}(z)~\mbox{and}~J_{f}(\xi)=\lim_{z\rightarrow
\xi\in\mathbb{T},z\in\mathbb{D}}J_{f}(z)$$ almost everywhere in
$\IT$ follows from  Theorem \Ref{KM-P}  and Lemma \ref{lem-main}.

For convenience, in the following, let
$$ C_{2}(K,\varphi,g)=\sup_{z\in\mathbb{D}}\|D_{f}(z)\|.$$

Since for almost all $z_1$ and $z_2\in \ID$,
\be\label{sun-7}
|f(z_1)-f(z_2)=\Big|\int_{[z_1,z_2]}f_zdz+f_{\overline{z}}d\overline{z} \Big|\leq C_{2}(K,\varphi,g) |z_1-z_2|,
\ee
we see that, to prove the Lipschitz continuity of $f$, it suffices to estimate the quantity $C_{2}(K,\varphi,g)$.
To reach this goal, we first show that the quantity $C_{2}(K,\varphi,g)$ satisfies an inequality which is stated in the following claim.
\bcl\label{eq-sh-18} $C_{2}(K,\varphi,g)\leq \big(C_{2}(K,\varphi,g)\big)^{1-\frac{1}{K}}\mu_1+\mu_2,$
 where $$\mu_1=\frac{K(Q(K))^{\frac{1}{K}+1}}{2\pi}\int_{0}^{2\pi}|1-e^{it}|^{-1+\frac{1}{K^{2}}}
dt,$$ $Q(K)$ is from Theorem \Ref{Mori}, $\mu_2=\mu_3+\mu_4,$
 $$\mu_3=\frac{K\|\varphi\|_{\infty}}{2}\left(\frac{\pi^{2}}{3}-1\right)^{\frac{1}{2}}
 +\frac{K\|g\|_{\infty}}{16}\left[1+2^{\frac{1}{2}}\left(1+\frac{\pi^{2}}{6}\right)^{\frac{1}{2}}\right],$$
 and
 $$\mu_4=\frac{\|\varphi\|_{\infty}}{2}\left[\max_{x\in[0,1)}\{h(x)\}+
2\left(\frac{\pi^{2}}{3}-1\right)^{\frac{1}{2}}\right]+\|g\|_{\infty}\left[\frac{53}{240}+\frac{2^{\frac{1}{2}}\left(1+\frac{\pi^{2}}{6}\right)^{\frac{1}{2}}}{8}\right].$$
\ecl

To prove the claim, we need the following preparation. Firstly, we
prove that for almost every $\theta\in [0,2\pi]$, \be\label{sun-3}
\|D_{f}(e^{i\theta})\|\leq\frac{K}{2\pi}\int_{0}^{2\pi}\frac{|f(e^{it})-f(e^{i\theta})|^{2}}{|e^{it}-e^{i\theta}|^{2}}dt+\mu_3.\ee

Since $f$ is a $K$-quasiconformal self-mapping of $\mathbb{D}$, we
see that $f$ can  be extended to the homeomorphism of
$\overline{\mathbb{D}}$ onto itself. For $\theta\in [0,2\pi]$, let
$$f(e^{i\theta})=f^{\ast}(e^{i\theta})=e^{i\eta(\theta)}.$$
Then, by (\ref{eq-d1}), we see that
$f(e^{i\theta})$ is absolutely continuous. It follows that
$$i\eta'(\theta)e^{i\eta(\theta)}=\frac{d}{d\theta}f(e^{i\theta})=\lim_{r\rightarrow1^{-}}\frac{\partial}{\partial\theta}f(re^{i\theta})=\lim_{r\rightarrow1^{-}}
\big[ir\left(f_{z}(re^{i\theta})e^{i\theta}-f_{\overline{z}}(re^{i\theta})e^{-i\theta}\right)\big],$$
which implies
\be\label{eq-x1.2}\frac{1}{K}\|D_{f}(e^{i\theta})\|\leq
\lim_{r\rightarrow1^{-}}\lambda(D_{f}(re^{i\theta}))\leq\eta'(\theta)\leq\lim_{r\rightarrow1^{-}}\|D_{f}(re^{i\theta})\|=\|D_{f}(e^{i\theta})\|\ee
almost everywhere in $[0,2\pi],$ where $r\in[0,1)$.

Since the existence of the two limits
$D_{f}(e^{i\theta})=\lim_{z\rightarrow
e^{i\theta},z\in\mathbb{D}}D_{f}(z)$ and $J_{f}(e^{i\theta})=\lim_{z\rightarrow
e^{i\theta},z\in\mathbb{D}}J_{f}(z)$ almost everywhere in
$[0,2\pi]$ guarantees that
$$\|D_{f}(e^{i\theta})\|^{2}\leq K J_{f}(e^{i\theta}),$$
we deduce from (\ref{eq-sh-1}) and (\ref{eq-x1.2}) that
\beq\label{eq-sh-11}\nonumber
 \|D_{f}(e^{i\theta})\|^{2}\leq K
\|D_{f}(e^{i\theta})\|\Bigg\{\frac{1}{2\pi}\int_{0}^{2\pi}\frac{|f(e^{it})-f(e^{i\theta})|^{2}}{|e^{it}-e^{i\theta}|^{2}}dt+\frac{\mu_3}{K}\Bigg\},\eeq
from which the inequality \eqref{sun-3} follows.

Secondly, we show that for any $\epsilon>0$, there exists
$\theta_{\epsilon}\in [0,2\pi]$ such that \be\label{sun-5}
C_{2}(K,\varphi,g)\leq(1+\epsilon)\|D_{f}(e^{i\theta_{\epsilon}})\|+\mu_4.\ee

For the proof, let $t\in[0,2\pi]$, and let
$$H_{t}(z)=\frac{\partial}{\partial
z}\mathcal{P}_{f^{\ast}}(z)+e^{it}\overline{\frac{\partial}{\partial
\overline{z}}\mathcal{P}_{f^{\ast}}(z)}$$ in $\mathbb{D}$.

Since
$\mathcal{P}_{f^{\ast}}=f-G_{1}[\varphi]+G_{2}[g]$ is harmonic, we
see that $H_{t}$ is analytic in
$\mathbb{D}$, and thus,
$$|H_{t}(z)|\leq\mbox{esssup}_{\theta\in[0,2\pi]}|H_{t}(e^{i\theta})|\leq\mbox{esssup}_{\theta\in[0,2\pi]}\|D_{\mathcal{P}_{f^{\ast}}}(e^{i\theta})\|.$$
Then, the facts
$$\|D_{\mathcal{P}_{f^{\ast}}}(z)\|=\max_{t\in[0,2\pi]}|H_{t}(z)|$$
and
$$\|D_{\mathcal{P}_{f^{\ast}}}(z)\|=\left|\frac{\partial f}{\partial z}-\frac{\partial G_{1}[\varphi]}{\partial z}+\frac{\partial G_{2}[g]}{\partial
z}\right|+\left|\frac{\partial f}{\partial
\overline{z}}-\frac{\partial G_{1}[\varphi]}{\partial
\overline{z}}+\frac{\partial G_{2}[g]}{\partial \overline{z}}\right|
$$
 ensure

\begin{eqnarray*}
\|D_{\mathcal{P}_{f^{\ast}}}(z)\|&\leq&
\mbox{esssup}_{\theta\in[0,2\pi]}\|D_{f}(e^{i\theta})\|+
\mbox{esssup}_{\theta\in[0,2\pi]}\|D_{G_{1}[\varphi]}(e^{i\theta})\|\\
&&+\mbox{esssup}_{\theta\in[0,2\pi]}\|D_{G_{2}[g]}(e^{i\theta})\|,
\end{eqnarray*} which, together with Lemmas \ref{lem-1} and
\ref{lem-2}, guarantees that for all $z\in\mathbb{D},$
$$
\|D_{f}(z)\|\leq\|D_{\mathcal{P}_{f^{\ast}}}(z)\|+\|D_{G_{1}[\varphi]}(z)\|+\|D_{G_{2}[g]}(z)\|\leq
\mbox{esssup}_{\theta\in[0,2\pi]}\|D_{f}(e^{i\theta})\|+\mu_4,
$$ from which the inequality \eqref{sun-5}  follows.

Let
$$\nu=\frac{1}{2\pi}\int_{0}^{2\pi}\frac{|f(e^{it})-f(e^{i\theta_{\epsilon}})|^{2}}{|e^{it}-e^{i\theta_{\epsilon}}|^{2}}dt.$$
Finally, we need the following estimate of $\nu$:
\be\label{sun-6}\nu\leq
\frac{\big(C_{2}(K,\varphi,g)\big)^{1-\frac{1}{K}}(Q(K))^{\frac{1}{K}+1}}{2\pi}\int_{0}^{2\pi}|e^{it}-e^{i\theta_{\epsilon}}|^{-1+\frac{1}{K^{2}}}
dt. \ee

Since it follows from \eqref{sun-7} that for almost all
$\theta_{1}$, $\theta_{2}\in [0,2\pi],$
\be\label{sun-8}|f(e^{i\theta_{1}})-f(e^{i\theta_{2}})|\leq
C_{2}(K,\varphi,g)\left|e^{i\theta_{1}}-e^{i\theta_{2}}\right|,\ee
we infer that
$$ \nu \leq
\frac{\big(C_{2}(K,\varphi,g)\big)^{1-\frac{1}{K}}}{2\pi}\int_{0}^{2\pi}|e^{it}-e^{i\theta_{\epsilon}}|^{-1+\frac{1}{K^{2}}}
\frac{|f(e^{it})-f(e^{i\theta_{\epsilon}})|^{1+\frac{1}{K}}}{|e^{it}-e^{i\theta_{\epsilon}}|^{\frac{1}{K}+\frac{1}{K^{2}}}}dt,
$$ from which, together with   Theorem \Ref{Mori}, the inequality \eqref{sun-6} follows.\medskip

Now, we are ready to finish the proof of the claim. It follows from
 \eqref{sun-5} that
$$C_{2}(K,\varphi,g)\leq(1+\epsilon)\|D_{f}(e^{i\theta_{\epsilon}})\|+\mu_4,$$
and so, \eqref{sun-3} and \eqref{sun-6} give \be\label{eq-sh-17}
C_{2}(K,\varphi,g) \leq
\big(C_{2}(K,\varphi,g)\big)^{1-\frac{1}{K}}\mu_1(1+\epsilon)+\mu_3(1+\epsilon)+\mu_4.
 \ee

 Moreover, by
\cite[Lemma 1.6]{K3}, we know that
$$\int_{0}^{2\pi}|e^{it}-e^{i\theta_{\epsilon}}|^{-1+\frac{1}{K^{2}}}
dt<\infty,$$ which shows $\mu_1<\infty$.

By letting $\epsilon\rightarrow0^{+}$, we get from (\ref{eq-sh-17})
that
$$C_{2}(K,\varphi,g)\leq
\big(C_{2}(K,\varphi,g)\big)^{1-\frac{1}{K}}\mu_1+\mu_2,$$ as required.\medskip

The following is a lower bound for $C_{2}(K,\varphi,g)$.
\bcl\label{sun-1} $C_{2}(K,\varphi,g)\geq 1$.\ecl

 Since
$$\int_{0}^{2\pi}\eta'(\theta)d\theta=\eta(2\pi)-\eta(0)=2\pi,$$
we conclude that
$$\mbox{esssup}_{\theta\in[0,2\pi]}\lim_{t\rightarrow\theta}\left|\frac{f(e^{i\theta})-f(e^{it})}{e^{i\theta}-e^{it}}\right|=\mbox{esssup}_{\theta\in[0,2\pi]}\eta'(\theta)\geq1.$$
Then, it follows from (\ref{sun-8}) and the following fact
$$\mbox{esssup}_{\theta\in[0,2\pi]}\lim_{t\rightarrow\theta}\left|\frac{f(e^{i\theta})-f(e^{it})}{e^{i\theta}-e^{it}}\right| \leq \mbox{esssup}_{0\leq\theta\neq
 t\leq 2\pi}\left|\frac{f(e^{i\theta})-f(e^{it})}{e^{i\theta}-e^{it}}\right|
 $$  that
$$C_{2}(K,\varphi,g)\geq 1.$$ Hence, the claim is true. \medskip

An upper bound of $C_{2}(K,\varphi,g)$ is established in the
following claim. \bcl\label{sun-2} If $\frac{(K-1)}{K}\mu_1<1$, then
$$C_{2}(K,\varphi,g)\leq \mu_5,$$ where
$\mu_5=\frac{\frac{1}{K}\mu_1+\mu_2}{1-\mu_1\left(1-\frac{1}{K}\right)}.$
\ecl
 The proof of this claim easily follows from \cite[Lemma 2.9]{K2}.\medskip

Now, we are ready to finish the discussions in this step. By Claims \ref{eq-sh-18} and \ref{sun-1}, we obtain
$$1\leq C_{2}(K,\varphi,g)\leq \mu_6,$$ where $\mu_6=(\mu_1+\mu_2)^{K}.$

By letting
$$C_{3}=\begin{cases}
\displaystyle \mu_6, &\mbox{ if }\, \frac{(K-1)}{K}\mu_1\geq1,\\
\displaystyle \min\{\mu_5, \mu_6\}, &\mbox{ if }\,
\frac{(K-1)}{K}\mu_1<1,
\end{cases}$$ we infer that
\be\label{eq-cw1} 1<C_{2}(K,\varphi,g)\leq C_{3}.\ee Then, the
Lipschtz continuity of $f$ easily follows from these estimates of
$C_{2}(K,\varphi,g)$.

\bst\label{bst-2} Co-Lipschitz continuity. \est We begin the
discussions of this step with some preparation which consists of the following two claims.

\bcl\label{eq-sh-28}
$\lambda(D_{\mathcal{P}_{f^{\ast}}}(e^{i\theta}))\geq
\frac{\mu_7}{K^{2}}-\left(1+\frac{1}{K^{2}}\right)\mu_8$ almost
everywhere in $[0,2\pi]$,
 \noindent where
 \be\label{thur-2} \mu_7=\max\{\mu'_7,~\mu''_7\},\;\;
\mu'_7=(Q(K))^{-2K}\frac{1}{2\pi}\int_{0}^{2\pi}|e^{it}-e^{i\theta}|^{2K-2}dt,
\ee
$$\mu''_7=\frac{1}{2}-\frac{\|\varphi\|_{\infty}}{8}-\frac{3\|g\|_{\infty}}{128},$$ and
\be\label{thur-1}
\mu_8=\frac{\|\varphi\|_{\infty}}{2}\left(\frac{\pi^{2}}{3}-1\right)^{\frac{1}{2}}
+\frac{\|g\|_{\infty}}{16}\left[1+2^{\frac{1}{2}}\left(1+\frac{\pi^{2}}{6}\right)^{\frac{1}{2}}\right],
\ee
\ecl

By (\ref{eq-x1.2}), we have

 $$\frac{\eta'(\theta)}{K}\leq\frac{\|D_{f}(e^{i\theta})\|}{K}\leq\lambda(D_{f}(e^{i\theta}))
 \leq\lambda(D_{\mathcal{P}_{f^{\ast}}}(e^{i\theta}))+\|D_{G_{1}[\varphi]}(e^{i\theta})\|+\|D_{G_{2}[g]}(e^{i\theta})\|,$$
which, together with Lemmas \ref{lem-1} and \ref{lem-2}, implies

\be\label{eerr-1}\lambda(D_{\mathcal{P}_{f^{\ast}}}(e^{i\theta}))\geq
\frac{\eta'(\theta)}{K}
-\|D_{G_{1}[\varphi]}(e^{i\theta})\|-\|D_{G_{2}[g]}(e^{i\theta})\|\geq\frac{\eta'(\theta)}{K}
-\mu_8.\ee
Then, we know from (\ref{eerr-1})  that, to prove the claim, it
suffices to show that \be\label{mon-2} K\eta'(\theta)\geq \mu_7.\ee

Again, it follows from (\ref{eq-x1.2}) that
$$\frac{J_{f}(e^{i\theta})}{\eta'(\theta)}\leq\frac{J_{f}(e^{i\theta})}{\lambda(D_{f}(e^{i\theta}))}
\leq
K\lambda(D_{f}(e^{i\theta}))\leq K\eta'(\theta),$$ and thus,
(\ref{eq-sh-1.1}) gives
$$K\eta'(\theta)\geq\frac{1}{2\pi}\int_{0}^{2\pi}\frac{|f(e^{it})-f(e^{i\theta})|^{2}}{|e^{it}-e^{i\theta}|^{2}}dt
-\mu_8. $$ This implies that, to prove \eqref{mon-2}, we only need
to verify the validity of the following inequality:
\be\label{eq-sh-25}\frac{1}{2\pi}\int_{0}^{2\pi}\frac{|f(e^{it})-f(e^{i\theta})|^{2}}{|e^{it}-e^{i\theta}|^{2}}dt
\geq \mu_7. \ee

We now prove this inequality. On the one hand, since $f^{-1}$ is  a
$K$-quasiconformal mapping, it follows from Theorem \Ref{Mori} that
for any $z_{1},z_{2}\in\mathbb{D}$, \beqq\label{eq-sh-26}
(Q(K))^{-K}|z_{1}-z_{2}|^{K}\leq|f(z_{1})-f(z_{2})|,\eeqq
which implies \beq\label{eq-sh-27}
\frac{1}{2\pi}\int_{0}^{2\pi}\frac{|f(e^{it})-f(e^{i\theta})|^{2}}{|e^{it}-e^{i\theta}|^{2}}dt\geq
\mu'_7.\eeq

On the other hand, since $f(0)=0$ and
$x^{2}\log\frac{1}{x}-(1-x^{2})\leq 0$ for $x\in(0,1]$, we see from
$$\big|G_{2}[g](0)\big|\leq \frac{\|g\|_{\infty}}{8\pi}\int_{\mathbb{D}}\left||\zeta|^{2}\log\frac{1}{|\zeta|}-(1-|\zeta|^{2})\right|d\sigma(\zeta)
=\frac{3\|g\|_{\infty}}{64}
$$
that \be\label{ccww-1} |\mathcal{P}_{f^{\ast}}(0)|
\leq\big|G_{1}[\varphi](0)\big|+\big|G_{2}[g](0)\big|\leq
\frac{\|\varphi\|_{\infty}}{4}+\frac{3\|g\|_{\infty}}{64}. \ee Then,
we infer from (\ref{ccww-1}) and the following fact:
$$\frac{1}{2\pi}\int_{0}^{2\pi}\frac{|f(e^{it})-f(e^{i\theta})|^{2}}{|e^{it}-e^{i\theta}|^{2}}dt\geq
\frac{1}{4\pi}\int_{0}^{2\pi}\left[1-\mbox{Re}\big(f(e^{it})\overline{f(e^{i\theta})}\big)\right]dt\geq\frac{1-|\mathcal{P}_{f^{\ast}}(0)|}{2}$$
that \beq\label{ccww-2}
\frac{1}{2\pi}\int_{0}^{2\pi}\frac{|f(e^{it})-f(e^{i\theta})|^{2}}{|e^{it}-e^{i\theta}|^{2}}dt
\geq \mu''_7.\eeq

Obviously, the inequality \eqref{eq-sh-25} follows from
\eqref{eq-sh-27} and \eqref{ccww-2}, and so, the claim is proved.
\bcl\label{ccww-6} For $z\in\mathbb{D},$
$\lambda(D_{\mathcal{P}_{f^{\ast}}}(z))\geq
\frac{\mu_7}{K^{2}}-\left(1+\frac{1}{K^{2}}\right)\mu_8.$ \ecl

By the Choquet-Rad\'o-Kneser theorem (see \cite{Cho}), we see that
$\mathcal{P}_{f^{\ast}}$ is a  sense-preserving  harmonic
diffeomorphism of $\mathbb{D}$ onto itself. Then, by Corollary \ref{H-lem}, we can let 


$$p_{1}(z)=\frac{\overline{\frac{\partial}{\partial
\overline{z}}\mathcal{P}_{f^{\ast}}(z)}}{\frac{\partial}{\partial
z}\mathcal{P}_{f^{\ast}}(z)}~\mbox{and}~p_{2}(z)=\left(\frac{\mu_7}{K^{2}}-\mu_8\right)\frac{1}{\frac{\partial}{\partial
z}\mathcal{P}_{f^{\ast}}(z)}$$
in $\ID$, and
for $\tau\in[0,2\pi]$, let
$$q_{\tau}(z)=p_{1}(z)+e^{i\tau}p_{2}(z).$$
Then, by Corollary \ref{H-lem},
\be\label{cw-ql}\sup_{z\in\mathbb{D}}|q_{\tau}(z)|<\infty\ee for all
$\tau\in[0,2\pi]$.

 By Claim \ref{eq-sh-28}, we have
\be\label{cw-qll}|q_{\tau}(e^{i\theta})|\leq|p_{1}(e^{i\theta})|+|p_{2}(e^{i\theta})|=
\frac{\left|\frac{\partial}{\partial
\overline{z}}\mathcal{P}_{f^{\ast}}(e^{i\theta})\right|+\frac{\mu_7}{K^{2}}-\left(1+\frac{1}{K^{2}}\right)\mu_8}
{\left|\frac{\partial}{\partial
z}\mathcal{P}_{f^{\ast}}(e^{i\theta})\right|}\leq1\ee almost
everywhere in $[0,2\pi]$. Let
$$E=\{\theta\in[0,2\pi]:~\lim_{z\rightarrow
e^{i\theta}}q_{\tau}(z)~\mbox{exists}\}.$$ Then  measure of
$[0,2\pi]\backslash E$ is zero. Hence it follows from (\ref{cw-ql})
and (\ref{cw-qll}) that, for all $\tau\in[0,2\pi]$,
\beqq\label{ccww-7}|q_{\tau}(z)|\leq
\frac{1}{2\pi}\int_{E}P(z,e^{i\theta})|q_{\tau}(e^{i\theta})|d\theta\leq
1,\eeqq from which, together with the   arbitrariness of
$\tau\in[0,2\pi]$, the claim follows.

Now, we are ready to finish the proof of the co-Lipschitz continuity
of $f$. Since
$$\lambda(D_{f}(z))\geq\lambda(D_{\mathcal{P}_{f^{\ast}}}(z))-\big(\|D_{G_{1}[\varphi]}(z)\|+\|D_{G_{2}[g]}(z)\|\big),$$
we see from Claim \ref{ccww-6}, Lemmas \ref{lem-1} and \ref{lem-2}
that \beq\label{mon-3}\lambda(D_{f}(z))\geq C_{1}(K,\varphi,g), \eeq
where \beq\label{ccww-10}
C_{1}(K,\varphi,g)&=&\frac{\mu_7}{K^{2}}-\frac{\|\varphi\|_{\infty}}{2}\left[\max_{x\in[0,1)}\{h(x)\}+
\left(\frac{\pi^{2}}{3}-1\right)^{\frac{1}{2}}\right]\\ \nonumber
&&-\left(1+\frac{1}{K^{2}}\right)\mu_8-\|g\|_{\infty}\left[\frac{19}{120}+\frac{2^{\frac{1}{2}}\left(1+\frac{\pi^{2}}{6}\right)^{\frac{1}{2}}}{16}\right].
\eeq And, we know from  \eqref{thur-2} and \eqref{thur-1} that
$C_{1}(K,\varphi,g)>0$ for small enough $\|g\|_{\infty}$ and
$\|\varphi\|_{\infty}$. Since for any $z_{1},$ $z_{2}\in\mathbb{D}$,
$$|f(z_{1})-f(z_{2})|\geq \int_{[z_{1},z_{2}]}\lambda(D_{f}(z))|dz|\geq C_{1}(K,\varphi,g)|z_{1}-z_{2}|,$$
 we conclude that $f$ is co-Lipschitz continuous.

\bst\label{bst-3} Bounds of the Lipschitz continuity coefficients $C_{1}(K,\varphi,g)$ and
$C_{2}(K,\varphi,g)$.\est

The discussions of this step consists of the following two claims.

\bcl There are constants $M_{2}(K)$ and $N_{2}(K,\varphi,g)$ such that
\begin{enumerate} \item
$C_{2}(K,\varphi,g)\leq M_{2}(K)+N_{2}(K,\varphi,g);$
\item
$\lim_{K\rightarrow1}M_{2}(K)=1$, and
\item
$\lim_{\|\varphi\|_{\infty}\rightarrow0,\|g\|_{\infty}\rightarrow0}N_{2}(K,\varphi,g)=0.$
\end{enumerate}\ecl
 From   (\ref{eq-cw1}), we see
that
$$ 1\leq C_{2}(K,\varphi,g)\leq C_3,$$
where
$$ C_3= \begin{cases}
\displaystyle (\mu_1+\mu_2)^{K}, &\mbox{ if }\, \frac{(K-1)}{K}\mu_1\geq 1,\\
\displaystyle \min\left\{(\mu_1+\mu_2)^{K},
\frac{\frac{1}{K}\mu_1+\mu_2}{1-\mu_1\left(1-\frac{1}{K}\right)}\right\},
&\mbox{ if }\, \frac{(K-1)}{K}\mu_1<1.
\end{cases}$$
Then, we have
$$C_3=\begin{cases}
\displaystyle M_{2}^{'}(K)+N_{2}^{'}(K,\varphi,g), &\mbox{ if }\, \frac{(K-1)}{K}\mu_1\geq 1,\\
\displaystyle \min\big\{M_{2}^{'}(K)+N_{2}^{'}(K,\varphi,g),
M_{2}^{''}(K)+N_{2}^{''}(K,\varphi,g)\big\}, &\mbox{ if }\,
\frac{(K-1)}{K}\mu_1< 1,
\end{cases}$$
where $M_{2}^{'}(K)=\mu_1^{K}$,
$M_{2}^{''}(K)=\frac{\mu_1}{K-\mu_1(K-1)}$, $N_{2}^{'}(K,\varphi,g)=(\mu_1+\mu_2)^{K}-\mu_1^{K}$, and
$$N_{2}^{''}(K,\varphi,g)=\frac{\mu_2}{1-\mu_1\big(1-\frac{1}{K}\big)}.$$
Let
$$M_{2}(K)=\max\{M_{2}^{'}(K), M_{2}^{''}(K)\}$$ and
$$N_{2}(K,\varphi,g)=\max\{N_{2}^{'}(K,\varphi,g), N_{2}^{''}(K,\varphi,g)\}.$$

It follows from the facts
$$\lim_{K\rightarrow1}M_{2}(K)=1\;\;\mbox{and}\;\; \lim_{\|\varphi\|_{\infty}\rightarrow0,\|g\|_{\infty}\rightarrow0}N_{2}(K,\varphi,g)=0$$
that these two constants are what we need, and so, the claim is proved.

\bcl\label{claim-3.8} There are constants $M_{1}(K)$ and
$N_{1}(K,\varphi,g)$ such that
\begin{enumerate} \item
$C_{1}(K,\varphi,g)\geq M_{1}(K)-N_{1}(K,\varphi,g);$
\item
$\lim_{K\rightarrow1}M_{1}(K)=1$, and
\item
 $\lim_{\|\varphi\|_{\infty}\rightarrow0,\|g\|_{\infty}\rightarrow0}N_{1}(K,\varphi,g)=0.$
\end{enumerate}
\ecl
By (\ref{ccww-10}), we have
$$C_{1}(K,\varphi,g)\geq M_{1}(K)-N_{1}(K,\varphi,g),$$
where
$$M_{1}(K)=K^{-2}(Q(K))^{-2K}\frac{1}{2\pi}\int_{0}^{2\pi}|e^{it}-e^{i\theta}|^{2K-2}dt$$
and
\begin{eqnarray*}
N_{1}(K,\varphi,g)&=&\frac{\|\varphi\|_{\infty}}{2}\left[\max_{x\in[0,1)}\{h(x)\}+\Big(2+\frac{1}{K^{2}}\Big)\left(\frac{\pi^{2}}{3}-1\right)^{\frac{1}{2}}\right]
\\&&+\|g\|_{\infty}\left[\frac{1}{16K^{2}}+\frac{53}{240}+\frac{2^{\frac{1}{2}}\big(1+2K^{2}\big)\left(1+\frac{\pi^{2}}{6}\right)^{\frac{1}{2}}}{16K^{2}}\right].
\end{eqnarray*} 
The following facts
$$\lim_{K\rightarrow1}M_{1}(K)=1\;\;\mbox{and}\;\;\lim_{\|\varphi\|_{\infty}\rightarrow0,\|g\|_{\infty}\rightarrow0}N_{1}(K,\varphi,g)=0$$
show that these two constants are what we want, and thus, the claim is true.

Now, by the discussions of Steps \ref{bst-1} $\sim$ \ref{bst-3}, we see that the theorem is proved. \qed \medskip

As a direct consequence of Claim \ref{claim-3.8}, we have the
following result.
\bcor\label{tue-1} Under the assumptions of
Theorem \ref{thm-1.1}, if, further, $\|g\|_{\infty}\leq a_1(K)$ and
$\|\varphi\|_{\infty}\leq a_2(K)$, then $f$ is co-Lipschitz
continuous, and so, it is bi-Lipschitz continuous, where
$a_1(K)=\frac{60}{(25+61K^{2})46^{2(K-1)}}$ and
$a_2(K)=\frac{25}{(38+101K^{2})46^{2(K-1)}}.$ \ecor

\section{Two examples}\label{sec-4}
The purpose of this section is to construct two examples. The first example shows that the co-Lipschitz continuity of $f$ from Theorem \ref{thm-1.1}
is invalid for arbitrary $g$ and $\varphi$, and the second one illustrates the possibility of $f$ to be bi-Lipschitz continuous.

\begin{example}\label{wed-5}
 For $z\in\overline{\mathbb{D}}$, let
$$f(z)=\beta|z|^{\gamma}z,$$
 where
$\gamma$ and $\beta$ are constants with $\gamma>3$ and $|\beta|=1$.
Then, $f$ is a four times continuously differentiable and
$K$-quasiconformal self-mapping of $\mathbb{D}$ with $f(0)=0$ and
$K=1+\gamma$. Furthermore, for $z\in\mathbb{D}$,
$$g(z)=\Delta(\Delta
f(z))=\beta\gamma^{2}(\gamma^{2}-4)\frac{|z|^{\gamma-2}}{\overline{z}},$$
and for $\xi\in \mathbb{T}$,
$$\varphi(\xi)=\beta\gamma(2+\gamma)\xi \;\;\mbox{and}\;\; f^*(\xi)=\beta\xi.$$
Then $\|g\|_{\infty}=\gamma^{2}(\gamma^{2}-4)$ and
$\|\varphi\|_{\infty}=\gamma(2+\gamma)$. However, $f$ is not
co-Lipschitz continuous (i.e., it does not satisfy \eqref{mon-1})
because
$$\lambda(D_{f}(0))=|f_{z}(0)|-|f_{\overline{z}}(0)|=0.$$
\end{example}

\begin{example}\label{wed-6}
For $z\in\overline{\mathbb{D}}$, let
$$f(z)=z+\frac{1}{200}(|z|^{2}-|z|^{4}).$$
Then, we have
$$g(z)=\Delta(\Delta
f(z))=-\frac{8}{25}$$ in $\ID$ and
$$\varphi(\xi)=-\frac{3}{50},\;\; f^*(\xi)=\xi$$ in $\mathbb{T}$.

Moreover, it is not difficult to know that $f$ is a
$K$-quasiconformal self-mapping of $\mathbb{D}$ with
$$
K=\max_{z\in\overline{\mathbb{D}}}\left\{\frac{|1+\overline{z}(1-2|z|^{2})M|+|Mz(1-2|z|^{2})|}{|1+\overline{z}(1-2|z|^{2})M|-|Mz(1-2|z|^{2})|}\right\}=\frac{100}{99},$$
where $M=\frac{1}{200}$. Since elementary computations lead to
$$a_1(K)=\frac{60}{(25+61K^{2})46^{2(K-1)}}>0.63\;\;\mbox{and}\;\; a_2(K)=\frac{25}{(38+101K^{2})46^{2(K-1)}}>0.16,$$
we know that
$$\|g\|_{\infty}<a_1(K)\;\;\mbox{and}\;\;\|\varphi\|_{\infty}< a_2(K),$$ where $a_i(K)$ $(i=1,$ $2)$ are from Corollary \ref{tue-1}.
Now, it follows from Corollary \ref{tue-1} that $f$ is co-Lipschitz continuous, and so, it is bi-Lipschitz continuous.
\end{example}

\normalsize

\end{document}